\documentclass[12pt]{article}
\usepackage{amsmath,amsthm,amssymb, amstext, amsopn, amsxtra}
\usepackage{latexsym}
\usepackage{amsfonts}
\begin{document}
\begin{center}
\Large{\bf Group Rings that are Additively Generated by Idempotents and Units}

\vspace{2mm}
\large{Dinesh Khurana and Chanchal Kumar} 
\end{center}

\vspace{5mm}
\noindent{\bf Abstract:} {\footnotesize Let $R$ be an Abelian\footnote{A ring is said to be Abelian if its
all idempotents are central} exchange ring. We prove the following results:\\
1. $R\mathbb{Z}_2$ and $RS_3$ are clean rings. \\
2. If $G$ is a group of prime order $p$ and $p$ is in the Jacobson radical of $R$, then 
$RG$ is clean. 3. If identity in $R$ is a sum of two units and $G$ is a locally
finite group, then every element in $RG$ is a sum of two units.\\
4. For any locally finite group $G$, $RG$ has stable range one.}

\vspace{6mm}

 All rings in this note are associative with identity. An element of a ring is said to be clean if it is a sum of a unit and an idempotent. 
A ring $R$ is said to be clean if its every element is clean.
These rings were introduced by Nicholson in [N$_1$] as a class of examples of exchange rings. 
In [N$_1$, Proposition 1.8] Nicholson proved that an Abelian exchange ring is clean. 

\vspace{2mm}
This work is motivated by the paper [M] of McGovern where it is proved that for a commutative clean ring $R$, the group ring $R\mathbb{Z}_2$ is clean. We extend this result by proving that $R\mathbb{Z}_2$ 
is clean whenever $R$ is an Abelian exchange ring.  Moreover our proof is quite short.
We also prove that $RS_3$ is clean for any Abelian exchange ring $R$. 

\vspace{2mm}
Let $R$ be a commutative clean ring and $G$ be a finite group of prime order $p$ such that
$p$ is invertible in $R$. In [HN, Example 1], Han and Nicholson gave an example to show that the group ring $RG$ may not be clean. We prove that if $R$ is an Abelian exchange ring and $G$ is a group of prime order $p$ such that $p \in J(R)$, then $RG$ is clean. 

\vspace{2mm} 
A lot of people have studied rings in which every element is a sum of two units (see [KS] and its references). An obvious necessary condition for the identity element of a ring $R$ to be a sum of two units is that $R$ does not have a factor ring isomorphic to $\mathbb{Z}_2$. In [KS, Theorem] it is proved  that if $R$ is a right self-injective ring which has no factor ring isomorphic to $\mathbb{Z}_2$, then every element of $R$ is a sum of two units. We prove that every element
of a group ring $RG$, where $R$ is an Abelian exchange ring and $G$ is a locally finite group, is a sum of two units whenever $RG$ has no factor ring isomorphic to $\mathbb{Z}_2$.

\vspace{2mm}
A ring $R$ is said to have stable range one if for any $a$, $b \in R$, $aR + bR = R$ implies that $a + bx$ is a unit
for some $x\in R$. This property was defined by Hyman Bass who also proved that every semilocal ring has stable range one. Evans proved that if the endomorphism ring of a module $M_R$ has stable one,  then $M_R$ has cancellation property in the category of right $R$-modules. For these results and more information, we refer the reader to [L]. We prove that if $R$ is an Abelian exchange ring and $G$ is a locally finite group, then $RG$ has stable range one.

\vspace{3mm}
The following result contains the crucial idea of the note. \\[4mm]
{\bf Lemma 1.} {\em If an element $x \in RG$ is not clean (resp., is not a sum of two units), then there exists
an ideal $I \neq R$ of $R$ such that $R/I$ does not have any nontrivial central idempotent and 
$\overline{x} \in (R/I)G$ is not clean (resp., is not a sum of two units).}\\[3mm]
{\bf Proof.} Suppose $x \in RG$ is not clean. Let $\cal{C}$ be the family of all ideals $I$ of $R$ such that $\overline{x} = x + IG \in RG/IG$ is not clean. Clearly $R\not \in \cal{C}$ and  $(0) \in \cal{C}$.
If $\{I_i\}$ is a chain of ideals in $\cal{C}$, then $I = \cup_iI_i$ is also in $\cal{C}$. Because otherwise $\overline{x} \in RG/IG$ and thus $\overline{x} \in RG/I_iG$, for some $i$, is clean. So $\cal{C}$ is inductive and so by Zorn's lemma, $\cal{C}$ has a maximal member,
say $I$. If $R/I$ has a non-trivial central idempotents, then $R/I = I_1/I \times I_2/I$, for
some ideals $I_1$ and $I_2$ of $R$ properly containing $I$. So $RG/IG \cong I_1G/IG \times I_2G/IG$. Now $\overline{x} \in RG/I_2G \cong RG/IG/I_2G/IG \cong I_1G/IG$ is clean by the maximality of $I$. Similarly, $\overline{x} \in I_2G/IG$ is clean implying that $\overline{x} \in RG/IG$ is clean. This is a contradiction. The other part, when $x$ is not a sum of two units, can be proved similarly. \qed

\vspace{2mm}
As an Abelian exchange ring without any non-trivial central idempotents is local, we get\\[3mm] 
{\bf Corollary 2.} {\em If $R$ is an Abelian exchange ring and $x\in RG$ is not clean (resp., is not a sum of two units), then there exists an ideal $I\neq R$ of $R$ such that $R/I$ is local and $\overline{x} \in RG/IG$ is not clean (resp., is not a sum of two unit).}

\vspace{3mm}
The part 1 of the following Lemma is due to Nicholson [N, Theorem] and part 2 is due to Woods
[W, Lemma 6.1].\\[3mm]
{\bf Lemma 3.} {\em 1. If $R$ is local, $G$ is locally finite $p$ group and $p \in J(R)$, then $RG$ is local\\[2mm]
2. If $R$ is semiperfect, then $RS_3$ is semiperfect.}

\vspace{3mm}
We are now ready to prove our first main result.\\[3mm]
{\bf Theorem 4.}  {\em If $R$ is an Abelian exchange ring, then $R\mathbb{Z}_2$ and $RS_3$ are clean.}\\[3mm]
{\bf Proof.} Let $G=\{1, g\}$. If $x \in RG$ is not clean, then by Corollary 2 there exists an ideal $I \neq R$
of $R$ such that $R/I$ is local and $\overline{x} \in RG/IG$ is not clean. We will show that
$RG/IG \cong (R/I)G$ is clean. Indeed if $2$ is a unit in $S = R/I$, then $a + bg \to (a+b, a-b)$ is an isomorphism from $SG$ to $S \times S$
and so $SG$ is clean. If $2 \in J(S)$, then $SG$ is local, and thus clean, by Lemma 3.

Similarly, if we suppose that $RS_3$ is not clean, we get a local factor ring $R/I$ of $R$ such that $(R/I)S_3$ is not clean. But $(R/I)S_3$ is semiperfect by Lemma 3 and is thus
clean by [CY, Proposition 7]. \qed 

\vspace{3mm}
The following result gives a class of clean group rings.\\[3mm]
{\bf Theorem 5.} {\em Let $R$ be an Abelian exchange ring and $G$ be a group of prime order $p$. If $p\in J(R)$, then $RG$ is clean.}\\[3mm]
{\bf Proof.} If $x\in RG$ is not clean, then by Corollary 2, there exists a local factor
ring $R/I$ of $R$ for some ideal $I \neq R$ such that $\overline{x} \in RG/IG$ is not clean. If $p$ is a unit in $R/I$, then $pr - 1 \in I$ for some $r \in R$. But as $p \in J(R)$, $pr - 1$ is a unit in $R$. Thus $I = R$, a contradiction. Thus $p \in J(R/I)$ and so by Lemma 3, $(R/I)G \cong RG/IG$ is a local and hence a clean ring. \qed

\vspace{3mm}
The following result is due to Woods [W, Lemma 4.1]. \\[3mm]
{\bf Lemma 6.} {\em If $G$ is a locally finite group, then for any ring $R$, $J(R) \subseteq J(RG)$. In particular, $J(R)G \subseteq J(RG)$.}

\vspace{3mm}
We will need the following easy lemma to prove our next main result.\\[3mm]
{\bf Lemma 7.} {\em Let $S$ be a unital subring of $R$ such that $S$ does not have any factor ring isomorphic to $\mathbb{Z}_2$. Then $R$ also does not have any factor ring isomorphic to
$\mathbb{Z}_2$}\\[3mm]
{\bf Proof.} Suppose there exists an epimorphism $f:R \to \mathbb{Z}_2$. As $f(1) \neq 0$,
the restriction of $f$ to $S$ is also an epimorphism from $S \to \mathbb{Z}_2$ \qed

\vspace{3mm} Let $R$ be an Abelian exchange ring and $G$ be a locally finite group. The following interesting result, in particular, shows that if identity in $RG$ is a sum of two units, then every element in $RG$ is a sum of two units.\\[3mm]
{\bf Theorem 8.} {\em For an Abelian exchange ring $R$ and a locally finite group $G$, the following conditions are equivalent:\\[2mm]
(i) Identity in $RG$ is a sum of two units.\\[2mm]
(ii) Identity in $R$ is a sum of two units in $R$. \\[2mm]
(iii) $R$ does not have any factor ring isomorphic to $\mathbb{Z}_2$.\\[2mm]
(iv) Every element in $RG$ is a sum of two units.}\\[3mm]
{\em Proof.} The implications (i) implies (ii); (ii) implies (iii) and (iv) implies (i) are clear. So we only have to prove the implication (iii) implies (iv). Suppose that $R$  does not have any factor ring isomorphic to $\mathbb{Z}_2$. By Lemma 7, $RG$ also does not have any factor ring isomorphic to $Z_2$. Suppose, to the contrary, that some $x \in RG$ is not a sum of two
units. Then by Corollary 2, there exists an ideal $I \neq R$ of $R$ such that $R/I$ is a local ring and $\overline{x} \in RG/IG$ is not a sum of two units. We will show that every element in 
$(R/I)G\cong RG/IG$ is a sum of two units. Let $S = R/I$. By Lemma 6, $J(S)G \subseteq J(SG)$. Thus every
element in $SG$ is a sum of two units if and only if every element in $SG/J(S)G$ is a sum
of two units. Now  $SG/J(S)G \cong (S/J(S))G$ and $S/J(S)$ is a division ring. Let $D = S/J(S)$.
if $y \in DG$, then as $G$ is locally finite, there exists a finite subgroup $H$ of $G$ such that $y \in DH$. By  [C, Theorem 4],  $DH$ is a self-injective ring. If $DH$ has a  factor ring isomorphic to $\mathbb{Z}_2$, then by Lemma 7, $D\cong \mathbb{Z}_2$. But as $D \cong R/I$, so $R$  has a factor ring isomorphic to $\mathbb{Z}_2$, which is a contradiction. Thus $DH$ does not have a factor ring isomorphic to $\mathbb{Z}_2$. So it follows by [KS, Theorem], that every element of $DH$ is a sum of two units. \qed

\vspace{3mm}
 We now prove our last result.\\[3mm]
{\bf Theorem 9.} {\em Let $R$ be an Abelian exchange ring and $G$ be a locally finite group. Then $RG$ has stable range one.}\\[3mm]
{\bf Proof.}  Let $aRG + bRG = RG$. We have to show that there exists a unit of the form 
$a + by$ for some $y \in RG$. Suppose not. Let $\cal{C}$ be the family of ideals $I$ of $R$ such that $a+by + IG$ is not a unit in $RG/IG$ for any $y \in R$. It is easy to see that $\cal{C}$ is inductive and so, as in Corollary 2, there exists an ideal $I \neq R$ of $R$ such that $R/I$ is local and  $a+by + IG$ is not a unit in $RG/IG$ for any $y \in R$. We will show that $(R/I)G \cong RG/IG$ has stable range one. As $\overline{a}\overline{RG} + \overline{b}\overline{RG} = \overline{RG}$, where $\overline{RG} = RG/IG$,  this will give us the desired contradiction. Let $S = R/I$ and $D = S/J(S)$. As $J(S)G \subseteq J(SG)$, so $SG$ has stable range one if and only if $SG/J(S)G \cong (S/J(S))G = DG$ has stable range one. Let
$x,\, x_1,\, y,\, y_1 \in DG$ be such that $xx_1 + yy_1 = 1$. As $G$ is locally finite, we can find
a finite subgroup $H$ of $G$ such that $x,\, x_1,\, y,\, y_1 \in DH$. But as $D$ is a division ring and $H$ is finite, 
$DH$ is Artinian and so has stable range one (see [L, Corollary 2.10]). Thus there exists
$w \in DH$ such that $x + yw$ is a unit in $DH$ and hence in $DG$ also. \qed\\[3mm]
{\bf References}\\[3mm]
[C] I. G. Connell, {\em On the group ring}, Canad. J. Math {\bf 15} (1963), 650-685.\\[2mm]
[CY] V. P. Camillo and H. P. Yu, {\em Exchange rings, units and idempotents}, Comm. Algebra {\bf 22} (1994), 4737-4749.\\[2mm]
[HN] J. Han and W. K. Nicholson, {\em Extensions of clean rings}, Comm. Algebra {\bf 29} (2001), 2589-2595.\\[2mm]
[KS] D. Khurana and A. K. Srivastava, {\em Right self-injective rings in which every element is a sum of two units}, J. Algebra Appl. {\bf 6} (2007), 281-286.\\[2mm]
[L] T. Y. Lam, {\em A crash course on stable range, cancellation, substitution, and exchange}, J. Algebra Appl. {\bf 3}  (2004), 301--343.\\[2mm]
[M] W. W. McGovern, {\em A characterization of commutative clean rings}, 
 Int. J. Math. Game Theory Algebra  {\bf 15}  (2006), 403--413. \\[2mm]
[N] W. K. Nicholson, {\em Local group rings}, Canad. Math. Bull. {\bf 15} (1972), 137-138.\\[2mm]
[N$_1$] W. K. Nicholson, {\em Lifting idempotents and exchange rings}, Trans. Amer. math. Soc. {\bf 229} (1977), 269-278.\\[2mm]
[W] S. M. Woods, {\em Some results on semiperfect group rings}, Canad. J. Math {\bf 26}  (1974), 121--129.\\[3mm]
 Faculty of Mathematics \\
Indian Inst.~of Science Education \& Research \\
 Chandigarh 160\,019, India\\
{\tt dkhurana@iisermohali.ac.in}\\
{\tt chanchal@iisermohali.ac.in}

\end{document}